\documentclass[12pt,dvips]{amsart}
\usepackage[usenames]{color}
\usepackage{tikz}
\usepackage{amsmath}

\newtheorem{lemma}{Lemma}[section]
\newtheorem{theorem}{Theorem}[section]

\newtheorem{proposition}{Proposition}[section]

\DeclareMathOperator{\ASM}{ASM}

\author{Robin Langer}
\title[Lambda determinants]{Lambda Determinants}

\begin{document}
\definecolor{AliceBlue}{rgb}{0.94,0.97,1.00}
\definecolor{AntiqueWhite1}{rgb}{1.00,0.94,0.86}
\definecolor{AntiqueWhite2}{rgb}{0.93,0.87,0.80}
\definecolor{AntiqueWhite3}{rgb}{0.80,0.75,0.69}
\definecolor{AntiqueWhite4}{rgb}{0.55,0.51,0.47}
\definecolor{AntiqueWhite}{rgb}{0.98,0.92,0.84}
\definecolor{BlanchedAlmond}{rgb}{1.00,0.92,0.80}
\definecolor{BlueViolet}{rgb}{0.54,0.17,0.89}
\definecolor{CadetBlue1}{rgb}{0.60,0.96,1.00}
\definecolor{CadetBlue2}{rgb}{0.56,0.90,0.93}
\definecolor{CadetBlue3}{rgb}{0.48,0.77,0.80}
\definecolor{CadetBlue4}{rgb}{0.33,0.53,0.55}
\definecolor{CadetBlue}{rgb}{0.37,0.62,0.63}
\definecolor{CornflowerBlue}{rgb}{0.39,0.58,0.93}
\definecolor{DarkBlue}{rgb}{0.00,0.00,0.55}
\definecolor{DarkCyan}{rgb}{0.00,0.55,0.55}
\definecolor{DarkGoldenrod1}{rgb}{1.00,0.73,0.06}
\definecolor{DarkGoldenrod2}{rgb}{0.93,0.68,0.05}
\definecolor{DarkGoldenrod3}{rgb}{0.80,0.58,0.05}
\definecolor{DarkGoldenrod4}{rgb}{0.55,0.40,0.03}
\definecolor{DarkGoldenrod}{rgb}{0.72,0.53,0.04}
\definecolor{DarkGray}{rgb}{0.66,0.66,0.66}
\definecolor{DarkGreen}{rgb}{0.00,0.39,0.00}
\definecolor{DarkGrey}{rgb}{0.66,0.66,0.66}
\definecolor{DarkKhaki}{rgb}{0.74,0.72,0.42}
\definecolor{DarkMagenta}{rgb}{0.55,0.00,0.55}
\definecolor{DarkOliveGreen1}{rgb}{0.79,1.00,0.44}
\definecolor{DarkOliveGreen2}{rgb}{0.74,0.93,0.41}
\definecolor{DarkOliveGreen3}{rgb}{0.64,0.80,0.35}
\definecolor{DarkOliveGreen4}{rgb}{0.43,0.55,0.24}
\definecolor{DarkOliveGreen}{rgb}{0.33,0.42,0.18}
\definecolor{DarkOrange1}{rgb}{1.00,0.50,0.00}
\definecolor{DarkOrange2}{rgb}{0.93,0.46,0.00}
\definecolor{DarkOrange3}{rgb}{0.80,0.40,0.00}
\definecolor{DarkOrange4}{rgb}{0.55,0.27,0.00}
\definecolor{DarkOrange}{rgb}{1.00,0.55,0.00}
\definecolor{DarkOrchid1}{rgb}{0.75,0.24,1.00}
\definecolor{DarkOrchid2}{rgb}{0.70,0.23,0.93}
\definecolor{DarkOrchid3}{rgb}{0.60,0.20,0.80}
\definecolor{DarkOrchid4}{rgb}{0.41,0.13,0.55}
\definecolor{DarkOrchid}{rgb}{0.60,0.20,0.80}
\definecolor{DarkRed}{rgb}{0.55,0.00,0.00}
\definecolor{DarkSalmon}{rgb}{0.91,0.59,0.48}
\definecolor{DarkSeaGreen1}{rgb}{0.76,1.00,0.76}
\definecolor{DarkSeaGreen2}{rgb}{0.71,0.93,0.71}
\definecolor{DarkSeaGreen3}{rgb}{0.61,0.80,0.61}
\definecolor{DarkSeaGreen4}{rgb}{0.41,0.55,0.41}
\definecolor{DarkSeaGreen}{rgb}{0.56,0.74,0.56}
\definecolor{DarkSlateBlue}{rgb}{0.28,0.24,0.55}
\definecolor{DarkSlateGray1}{rgb}{0.59,1.00,1.00}
\definecolor{DarkSlateGray2}{rgb}{0.55,0.93,0.93}
\definecolor{DarkSlateGray3}{rgb}{0.47,0.80,0.80}
\definecolor{DarkSlateGray4}{rgb}{0.32,0.55,0.55}
\definecolor{DarkSlateGray}{rgb}{0.18,0.31,0.31}
\definecolor{DarkSlateGrey}{rgb}{0.18,0.31,0.31}
\definecolor{DarkTurquoise}{rgb}{0.00,0.81,0.82}
\definecolor{DarkViolet}{rgb}{0.58,0.00,0.83}
\definecolor{DeepPink1}{rgb}{1.00,0.08,0.58}
\definecolor{DeepPink2}{rgb}{0.93,0.07,0.54}
\definecolor{DeepPink3}{rgb}{0.80,0.06,0.46}
\definecolor{DeepPink4}{rgb}{0.55,0.04,0.31}
\definecolor{DeepPink}{rgb}{1.00,0.08,0.58}
\definecolor{DeepSkyBlue1}{rgb}{0.00,0.75,1.00}
\definecolor{DeepSkyBlue2}{rgb}{0.00,0.70,0.93}
\definecolor{DeepSkyBlue3}{rgb}{0.00,0.60,0.80}
\definecolor{DeepSkyBlue4}{rgb}{0.00,0.41,0.55}
\definecolor{DeepSkyBlue}{rgb}{0.00,0.75,1.00}
\definecolor{DimGray}{rgb}{0.41,0.41,0.41}
\definecolor{DimGrey}{rgb}{0.41,0.41,0.41}
\definecolor{DodgerBlue1}{rgb}{0.12,0.56,1.00}
\definecolor{DodgerBlue2}{rgb}{0.11,0.53,0.93}
\definecolor{DodgerBlue3}{rgb}{0.09,0.45,0.80}
\definecolor{DodgerBlue4}{rgb}{0.06,0.31,0.55}
\definecolor{DodgerBlue}{rgb}{0.12,0.56,1.00}
\definecolor{FloralWhite}{rgb}{1.00,0.98,0.94}
\definecolor{ForestGreen}{rgb}{0.13,0.55,0.13}
\definecolor{GhostWhite}{rgb}{0.97,0.97,1.00}
\definecolor{GreenYellow}{rgb}{0.68,1.00,0.18}
\definecolor{HotPink1}{rgb}{1.00,0.43,0.71}
\definecolor{HotPink2}{rgb}{0.93,0.42,0.65}
\definecolor{HotPink3}{rgb}{0.80,0.38,0.56}
\definecolor{HotPink4}{rgb}{0.55,0.23,0.38}
\definecolor{HotPink}{rgb}{1.00,0.41,0.71}
\definecolor{IndianRed1}{rgb}{1.00,0.42,0.42}
\definecolor{IndianRed2}{rgb}{0.93,0.39,0.39}
\definecolor{IndianRed3}{rgb}{0.80,0.33,0.33}
\definecolor{IndianRed4}{rgb}{0.55,0.23,0.23}
\definecolor{IndianRed}{rgb}{0.80,0.36,0.36}
\definecolor{LavenderBlush1}{rgb}{1.00,0.94,0.96}
\definecolor{LavenderBlush2}{rgb}{0.93,0.88,0.90}
\definecolor{LavenderBlush3}{rgb}{0.80,0.76,0.77}
\definecolor{LavenderBlush4}{rgb}{0.55,0.51,0.53}
\definecolor{LavenderBlush}{rgb}{1.00,0.94,0.96}
\definecolor{LawnGreen}{rgb}{0.49,0.99,0.00}
\definecolor{LemonChiffon1}{rgb}{1.00,0.98,0.80}
\definecolor{LemonChiffon2}{rgb}{0.93,0.91,0.75}
\definecolor{LemonChiffon3}{rgb}{0.80,0.79,0.65}
\definecolor{LemonChiffon4}{rgb}{0.55,0.54,0.44}
\definecolor{LemonChiffon}{rgb}{1.00,0.98,0.80}
\definecolor{LightBlue1}{rgb}{0.75,0.94,1.00}
\definecolor{LightBlue2}{rgb}{0.70,0.87,0.93}
\definecolor{LightBlue3}{rgb}{0.60,0.75,0.80}
\definecolor{LightBlue4}{rgb}{0.41,0.51,0.55}
\definecolor{LightBlue}{rgb}{0.68,0.85,0.90}
\definecolor{LightCoral}{rgb}{0.94,0.50,0.50}
\definecolor{LightCyan1}{rgb}{0.88,1.00,1.00}
\definecolor{LightCyan2}{rgb}{0.82,0.93,0.93}
\definecolor{LightCyan3}{rgb}{0.71,0.80,0.80}
\definecolor{LightCyan4}{rgb}{0.48,0.55,0.55}
\definecolor{LightCyan}{rgb}{0.88,1.00,1.00}
\definecolor{LightGoldenrod1}{rgb}{1.00,0.93,0.55}
\definecolor{LightGoldenrod2}{rgb}{0.93,0.86,0.51}
\definecolor{LightGoldenrod3}{rgb}{0.80,0.75,0.44}
\definecolor{LightGoldenrod4}{rgb}{0.55,0.51,0.30}
\definecolor{LightGoldenrodYellow}{rgb}{0.98,0.98,0.82}
\definecolor{LightGoldenrod}{rgb}{0.93,0.87,0.51}
\definecolor{LightGray}{rgb}{0.83,0.83,0.83}
\definecolor{LightGreen}{rgb}{0.56,0.93,0.56}
\definecolor{LightGrey}{rgb}{0.83,0.83,0.83}
\definecolor{LightPink1}{rgb}{1.00,0.68,0.73}
\definecolor{LightPink2}{rgb}{0.93,0.64,0.68}
\definecolor{LightPink3}{rgb}{0.80,0.55,0.58}
\definecolor{LightPink4}{rgb}{0.55,0.37,0.40}
\definecolor{LightPink}{rgb}{1.00,0.71,0.76}
\definecolor{LightSalmon1}{rgb}{1.00,0.63,0.48}
\definecolor{LightSalmon2}{rgb}{0.93,0.58,0.45}
\definecolor{LightSalmon3}{rgb}{0.80,0.51,0.38}
\definecolor{LightSalmon4}{rgb}{0.55,0.34,0.26}
\definecolor{LightSalmon}{rgb}{1.00,0.63,0.48}
\definecolor{LightSeaGreen}{rgb}{0.13,0.70,0.67}
\definecolor{LightSkyBlue1}{rgb}{0.69,0.89,1.00}
\definecolor{LightSkyBlue2}{rgb}{0.64,0.83,0.93}
\definecolor{LightSkyBlue3}{rgb}{0.55,0.71,0.80}
\definecolor{LightSkyBlue4}{rgb}{0.38,0.48,0.55}
\definecolor{LightSkyBlue}{rgb}{0.53,0.81,0.98}
\definecolor{LightSlateBlue}{rgb}{0.52,0.44,1.00}
\definecolor{LightSlateGray}{rgb}{0.47,0.53,0.60}
\definecolor{LightSlateGrey}{rgb}{0.47,0.53,0.60}
\definecolor{LightSteelBlue1}{rgb}{0.79,0.88,1.00}
\definecolor{LightSteelBlue2}{rgb}{0.74,0.82,0.93}
\definecolor{LightSteelBlue3}{rgb}{0.64,0.71,0.80}
\definecolor{LightSteelBlue4}{rgb}{0.43,0.48,0.55}
\definecolor{LightSteelBlue}{rgb}{0.69,0.77,0.87}
\definecolor{LightYellow1}{rgb}{1.00,1.00,0.88}
\definecolor{LightYellow2}{rgb}{0.93,0.93,0.82}
\definecolor{LightYellow3}{rgb}{0.80,0.80,0.71}
\definecolor{LightYellow4}{rgb}{0.55,0.55,0.48}
\definecolor{LightYellow}{rgb}{1.00,1.00,0.88}
\definecolor{LimeGreen}{rgb}{0.20,0.80,0.20}
\definecolor{MediumAquamarine}{rgb}{0.40,0.80,0.67}
\definecolor{MediumBlue}{rgb}{0.00,0.00,0.80}
\definecolor{MediumOrchid1}{rgb}{0.88,0.40,1.00}
\definecolor{MediumOrchid2}{rgb}{0.82,0.37,0.93}
\definecolor{MediumOrchid3}{rgb}{0.71,0.32,0.80}
\definecolor{MediumOrchid4}{rgb}{0.48,0.22,0.55}
\definecolor{MediumOrchid}{rgb}{0.73,0.33,0.83}
\definecolor{MediumPurple1}{rgb}{0.67,0.51,1.00}
\definecolor{MediumPurple2}{rgb}{0.62,0.47,0.93}
\definecolor{MediumPurple3}{rgb}{0.54,0.41,0.80}
\definecolor{MediumPurple4}{rgb}{0.36,0.28,0.55}
\definecolor{MediumPurple}{rgb}{0.58,0.44,0.86}
\definecolor{MediumSeaGreen}{rgb}{0.24,0.70,0.44}
\definecolor{MediumSlateBlue}{rgb}{0.48,0.41,0.93}
\definecolor{MediumSpringGreen}{rgb}{0.00,0.98,0.60}
\definecolor{MediumTurquoise}{rgb}{0.28,0.82,0.80}
\definecolor{MediumVioletRed}{rgb}{0.78,0.08,0.52}
\definecolor{MidnightBlue}{rgb}{0.10,0.10,0.44}
\definecolor{MintCream}{rgb}{0.96,1.00,0.98}
\definecolor{MistyRose1}{rgb}{1.00,0.89,0.88}
\definecolor{MistyRose2}{rgb}{0.93,0.84,0.82}
\definecolor{MistyRose3}{rgb}{0.80,0.72,0.71}
\definecolor{MistyRose4}{rgb}{0.55,0.49,0.48}
\definecolor{MistyRose}{rgb}{1.00,0.89,0.88}
\definecolor{NavajoWhite1}{rgb}{1.00,0.87,0.68}
\definecolor{NavajoWhite2}{rgb}{0.93,0.81,0.63}
\definecolor{NavajoWhite3}{rgb}{0.80,0.70,0.55}
\definecolor{NavajoWhite4}{rgb}{0.55,0.47,0.37}
\definecolor{NavajoWhite}{rgb}{1.00,0.87,0.68}
\definecolor{NavyBlue}{rgb}{0.00,0.00,0.50}
\definecolor{OldLace}{rgb}{0.99,0.96,0.90}
\definecolor{OliveDrab1}{rgb}{0.75,1.00,0.24}
\definecolor{OliveDrab2}{rgb}{0.70,0.93,0.23}
\definecolor{OliveDrab3}{rgb}{0.60,0.80,0.20}
\definecolor{OliveDrab4}{rgb}{0.41,0.55,0.13}
\definecolor{OliveDrab}{rgb}{0.42,0.56,0.14}
\definecolor{OrangeRed1}{rgb}{1.00,0.27,0.00}
\definecolor{OrangeRed2}{rgb}{0.93,0.25,0.00}
\definecolor{OrangeRed3}{rgb}{0.80,0.22,0.00}
\definecolor{OrangeRed4}{rgb}{0.55,0.15,0.00}
\definecolor{OrangeRed}{rgb}{1.00,0.27,0.00}
\definecolor{PaleGoldenrod}{rgb}{0.93,0.91,0.67}
\definecolor{PaleGreen1}{rgb}{0.60,1.00,0.60}
\definecolor{PaleGreen2}{rgb}{0.56,0.93,0.56}
\definecolor{PaleGreen3}{rgb}{0.49,0.80,0.49}
\definecolor{PaleGreen4}{rgb}{0.33,0.55,0.33}
\definecolor{PaleGreen}{rgb}{0.60,0.98,0.60}
\definecolor{PaleTurquoise1}{rgb}{0.73,1.00,1.00}
\definecolor{PaleTurquoise2}{rgb}{0.68,0.93,0.93}
\definecolor{PaleTurquoise3}{rgb}{0.59,0.80,0.80}
\definecolor{PaleTurquoise4}{rgb}{0.40,0.55,0.55}
\definecolor{PaleTurquoise}{rgb}{0.69,0.93,0.93}
\definecolor{PaleVioletRed1}{rgb}{1.00,0.51,0.67}
\definecolor{PaleVioletRed2}{rgb}{0.93,0.47,0.62}
\definecolor{PaleVioletRed3}{rgb}{0.80,0.41,0.54}
\definecolor{PaleVioletRed4}{rgb}{0.55,0.28,0.36}
\definecolor{PaleVioletRed}{rgb}{0.86,0.44,0.58}
\definecolor{PapayaWhip}{rgb}{1.00,0.94,0.84}
\definecolor{PeachPuff1}{rgb}{1.00,0.85,0.73}
\definecolor{PeachPuff2}{rgb}{0.93,0.80,0.68}
\definecolor{PeachPuff3}{rgb}{0.80,0.69,0.58}
\definecolor{PeachPuff4}{rgb}{0.55,0.47,0.40}
\definecolor{PeachPuff}{rgb}{1.00,0.85,0.73}
\definecolor{PowderBlue}{rgb}{0.69,0.88,0.90}
\definecolor{RosyBrown1}{rgb}{1.00,0.76,0.76}
\definecolor{RosyBrown2}{rgb}{0.93,0.71,0.71}
\definecolor{RosyBrown3}{rgb}{0.80,0.61,0.61}
\definecolor{RosyBrown4}{rgb}{0.55,0.41,0.41}
\definecolor{RosyBrown}{rgb}{0.74,0.56,0.56}
\definecolor{RoyalBlue1}{rgb}{0.28,0.46,1.00}
\definecolor{RoyalBlue2}{rgb}{0.26,0.43,0.93}
\definecolor{RoyalBlue3}{rgb}{0.23,0.37,0.80}
\definecolor{RoyalBlue4}{rgb}{0.15,0.25,0.55}
\definecolor{RoyalBlue}{rgb}{0.25,0.41,0.88}
\definecolor{SaddleBrown}{rgb}{0.55,0.27,0.07}
\definecolor{SandyBrown}{rgb}{0.96,0.64,0.38}
\definecolor{SeaGreen1}{rgb}{0.33,1.00,0.62}
\definecolor{SeaGreen2}{rgb}{0.31,0.93,0.58}
\definecolor{SeaGreen3}{rgb}{0.26,0.80,0.50}
\definecolor{SeaGreen4}{rgb}{0.18,0.55,0.34}
\definecolor{SeaGreen}{rgb}{0.18,0.55,0.34}
\definecolor{SkyBlue1}{rgb}{0.53,0.81,1.00}
\definecolor{SkyBlue2}{rgb}{0.49,0.75,0.93}
\definecolor{SkyBlue3}{rgb}{0.42,0.65,0.80}
\definecolor{SkyBlue4}{rgb}{0.29,0.44,0.55}
\definecolor{SkyBlue}{rgb}{0.53,0.81,0.92}
\definecolor{SlateBlue1}{rgb}{0.51,0.44,1.00}
\definecolor{SlateBlue2}{rgb}{0.48,0.40,0.93}
\definecolor{SlateBlue3}{rgb}{0.41,0.35,0.80}
\definecolor{SlateBlue4}{rgb}{0.28,0.24,0.55}
\definecolor{SlateBlue}{rgb}{0.42,0.35,0.80}
\definecolor{SlateGray1}{rgb}{0.78,0.89,1.00}
\definecolor{SlateGray2}{rgb}{0.73,0.83,0.93}
\definecolor{SlateGray3}{rgb}{0.62,0.71,0.80}
\definecolor{SlateGray4}{rgb}{0.42,0.48,0.55}
\definecolor{SlateGray}{rgb}{0.44,0.50,0.56}
\definecolor{SlateGrey}{rgb}{0.44,0.50,0.56}
\definecolor{SpringGreen1}{rgb}{0.00,1.00,0.50}
\definecolor{SpringGreen2}{rgb}{0.00,0.93,0.46}
\definecolor{SpringGreen3}{rgb}{0.00,0.80,0.40}
\definecolor{SpringGreen4}{rgb}{0.00,0.55,0.27}
\definecolor{SpringGreen}{rgb}{0.00,1.00,0.50}
\definecolor{SteelBlue1}{rgb}{0.39,0.72,1.00}
\definecolor{SteelBlue2}{rgb}{0.36,0.67,0.93}
\definecolor{SteelBlue3}{rgb}{0.31,0.58,0.80}
\definecolor{SteelBlue4}{rgb}{0.21,0.39,0.55}
\definecolor{SteelBlue}{rgb}{0.27,0.51,0.71}
\definecolor{VioletRed1}{rgb}{1.00,0.24,0.59}
\definecolor{VioletRed2}{rgb}{0.93,0.23,0.55}
\definecolor{VioletRed3}{rgb}{0.80,0.20,0.47}
\definecolor{VioletRed4}{rgb}{0.55,0.13,0.32}
\definecolor{VioletRed}{rgb}{0.82,0.13,0.56}
\definecolor{WhiteSmoke}{rgb}{0.96,0.96,0.96}
\definecolor{YellowGreen}{rgb}{0.60,0.80,0.20}
\definecolor{aliceblue}{rgb}{0.94,0.97,1.00}
\definecolor{antiquewhite}{rgb}{0.98,0.92,0.84}
\definecolor{aquamarine1}{rgb}{0.50,1.00,0.83}
\definecolor{aquamarine2}{rgb}{0.46,0.93,0.78}
\definecolor{aquamarine3}{rgb}{0.40,0.80,0.67}
\definecolor{aquamarine4}{rgb}{0.27,0.55,0.45}
\definecolor{aquamarine}{rgb}{0.50,1.00,0.83}
\definecolor{azure1}{rgb}{0.94,1.00,1.00}
\definecolor{azure2}{rgb}{0.88,0.93,0.93}
\definecolor{azure3}{rgb}{0.76,0.80,0.80}
\definecolor{azure4}{rgb}{0.51,0.55,0.55}
\definecolor{azure}{rgb}{0.94,1.00,1.00}
\definecolor{beige}{rgb}{0.96,0.96,0.86}
\definecolor{bisque1}{rgb}{1.00,0.89,0.77}
\definecolor{bisque2}{rgb}{0.93,0.84,0.72}
\definecolor{bisque3}{rgb}{0.80,0.72,0.62}
\definecolor{bisque4}{rgb}{0.55,0.49,0.42}
\definecolor{bisque}{rgb}{1.00,0.89,0.77}
\definecolor{black}{rgb}{0.00,0.00,0.00}
\definecolor{blanchedalmond}{rgb}{1.00,0.92,0.80}
\definecolor{blue1}{rgb}{0.00,0.00,1.00}
\definecolor{blue2}{rgb}{0.00,0.00,0.93}
\definecolor{blue3}{rgb}{0.00,0.00,0.80}
\definecolor{blue4}{rgb}{0.00,0.00,0.55}
\definecolor{blueviolet}{rgb}{0.54,0.17,0.89}
\definecolor{blue}{rgb}{0.00,0.00,1.00}
\definecolor{brown1}{rgb}{1.00,0.25,0.25}
\definecolor{brown2}{rgb}{0.93,0.23,0.23}
\definecolor{brown3}{rgb}{0.80,0.20,0.20}
\definecolor{brown4}{rgb}{0.55,0.14,0.14}
\definecolor{brown}{rgb}{0.65,0.16,0.16}
\definecolor{burlywood1}{rgb}{1.00,0.83,0.61}
\definecolor{burlywood2}{rgb}{0.93,0.77,0.57}
\definecolor{burlywood3}{rgb}{0.80,0.67,0.49}
\definecolor{burlywood4}{rgb}{0.55,0.45,0.33}
\definecolor{burlywood}{rgb}{0.87,0.72,0.53}
\definecolor{cadetblue}{rgb}{0.37,0.62,0.63}
\definecolor{chartreuse1}{rgb}{0.50,1.00,0.00}
\definecolor{chartreuse2}{rgb}{0.46,0.93,0.00}
\definecolor{chartreuse3}{rgb}{0.40,0.80,0.00}
\definecolor{chartreuse4}{rgb}{0.27,0.55,0.00}
\definecolor{chartreuse}{rgb}{0.50,1.00,0.00}
\definecolor{chocolate1}{rgb}{1.00,0.50,0.14}
\definecolor{chocolate2}{rgb}{0.93,0.46,0.13}
\definecolor{chocolate3}{rgb}{0.80,0.40,0.11}
\definecolor{chocolate4}{rgb}{0.55,0.27,0.07}
\definecolor{chocolate}{rgb}{0.82,0.41,0.12}
\definecolor{coral1}{rgb}{1.00,0.45,0.34}
\definecolor{coral2}{rgb}{0.93,0.42,0.31}
\definecolor{coral3}{rgb}{0.80,0.36,0.27}
\definecolor{coral4}{rgb}{0.55,0.24,0.18}
\definecolor{coral}{rgb}{1.00,0.50,0.31}
\definecolor{cornflowerblue}{rgb}{0.39,0.58,0.93}
\definecolor{cornsilk1}{rgb}{1.00,0.97,0.86}
\definecolor{cornsilk2}{rgb}{0.93,0.91,0.80}
\definecolor{cornsilk3}{rgb}{0.80,0.78,0.69}
\definecolor{cornsilk4}{rgb}{0.55,0.53,0.47}
\definecolor{cornsilk}{rgb}{1.00,0.97,0.86}
\definecolor{cyan1}{rgb}{0.00,1.00,1.00}
\definecolor{cyan2}{rgb}{0.00,0.93,0.93}
\definecolor{cyan3}{rgb}{0.00,0.80,0.80}
\definecolor{cyan4}{rgb}{0.00,0.55,0.55}
\definecolor{cyan}{rgb}{0.00,1.00,1.00}
\definecolor{darkblue}{rgb}{0.00,0.00,0.55}
\definecolor{darkcyan}{rgb}{0.00,0.55,0.55}
\definecolor{darkgoldenrod}{rgb}{0.72,0.53,0.04}
\definecolor{darkgray}{rgb}{0.66,0.66,0.66}
\definecolor{darkgreen}{rgb}{0.00,0.39,0.00}
\definecolor{darkgrey}{rgb}{0.66,0.66,0.66}
\definecolor{darkkhaki}{rgb}{0.74,0.72,0.42}
\definecolor{darkmagenta}{rgb}{0.55,0.00,0.55}
\definecolor{darkolive}{rgb}{0.33,0.42,0.18}
\definecolor{darkorange}{rgb}{1.00,0.55,0.00}
\definecolor{darkorchid}{rgb}{0.60,0.20,0.80}
\definecolor{darkred}{rgb}{0.55,0.00,0.00}
\definecolor{darksalmon}{rgb}{0.91,0.59,0.48}
\definecolor{darksea}{rgb}{0.56,0.74,0.56}
\definecolor{darkslate}{rgb}{0.18,0.31,0.31}
\definecolor{darkslate}{rgb}{0.18,0.31,0.31}
\definecolor{darkslate}{rgb}{0.28,0.24,0.55}
\definecolor{darkturquoise}{rgb}{0.00,0.81,0.82}
\definecolor{darkviolet}{rgb}{0.58,0.00,0.83}
\definecolor{deeppink}{rgb}{1.00,0.08,0.58}
\definecolor{deepsky}{rgb}{0.00,0.75,1.00}
\definecolor{dimgray}{rgb}{0.41,0.41,0.41}
\definecolor{dimgrey}{rgb}{0.41,0.41,0.41}
\definecolor{dodgerblue}{rgb}{0.12,0.56,1.00}
\definecolor{firebrick1}{rgb}{1.00,0.19,0.19}
\definecolor{firebrick2}{rgb}{0.93,0.17,0.17}
\definecolor{firebrick3}{rgb}{0.80,0.15,0.15}
\definecolor{firebrick4}{rgb}{0.55,0.10,0.10}
\definecolor{firebrick}{rgb}{0.70,0.13,0.13}
\definecolor{floralwhite}{rgb}{1.00,0.98,0.94}
\definecolor{forestgreen}{rgb}{0.13,0.55,0.13}
\definecolor{gainsboro}{rgb}{0.86,0.86,0.86}
\definecolor{ghostwhite}{rgb}{0.97,0.97,1.00}
\definecolor{gold1}{rgb}{1.00,0.84,0.00}
\definecolor{gold2}{rgb}{0.93,0.79,0.00}
\definecolor{gold3}{rgb}{0.80,0.68,0.00}
\definecolor{gold4}{rgb}{0.55,0.46,0.00}
\definecolor{goldenrod1}{rgb}{1.00,0.76,0.15}
\definecolor{goldenrod2}{rgb}{0.93,0.71,0.13}
\definecolor{goldenrod3}{rgb}{0.80,0.61,0.11}
\definecolor{goldenrod4}{rgb}{0.55,0.41,0.08}
\definecolor{goldenrod}{rgb}{0.85,0.65,0.13}
\definecolor{gold}{rgb}{1.00,0.84,0.00}
\definecolor{gray0}{rgb}{0.00,0.00,0.00}
\definecolor{gray100}{rgb}{1.00,1.00,1.00}
\definecolor{gray10}{rgb}{0.10,0.10,0.10}
\definecolor{gray11}{rgb}{0.11,0.11,0.11}
\definecolor{gray12}{rgb}{0.12,0.12,0.12}
\definecolor{gray13}{rgb}{0.13,0.13,0.13}
\definecolor{gray14}{rgb}{0.14,0.14,0.14}
\definecolor{gray15}{rgb}{0.15,0.15,0.15}
\definecolor{gray16}{rgb}{0.16,0.16,0.16}
\definecolor{gray17}{rgb}{0.17,0.17,0.17}
\definecolor{gray18}{rgb}{0.18,0.18,0.18}
\definecolor{gray19}{rgb}{0.19,0.19,0.19}
\definecolor{gray1}{rgb}{0.01,0.01,0.01}
\definecolor{gray20}{rgb}{0.20,0.20,0.20}
\definecolor{gray21}{rgb}{0.21,0.21,0.21}
\definecolor{gray22}{rgb}{0.22,0.22,0.22}
\definecolor{gray23}{rgb}{0.23,0.23,0.23}
\definecolor{gray24}{rgb}{0.24,0.24,0.24}
\definecolor{gray25}{rgb}{0.25,0.25,0.25}
\definecolor{gray26}{rgb}{0.26,0.26,0.26}
\definecolor{gray27}{rgb}{0.27,0.27,0.27}
\definecolor{gray28}{rgb}{0.28,0.28,0.28}
\definecolor{gray29}{rgb}{0.29,0.29,0.29}
\definecolor{gray2}{rgb}{0.02,0.02,0.02}
\definecolor{gray30}{rgb}{0.30,0.30,0.30}
\definecolor{gray31}{rgb}{0.31,0.31,0.31}
\definecolor{gray32}{rgb}{0.32,0.32,0.32}
\definecolor{gray33}{rgb}{0.33,0.33,0.33}
\definecolor{gray34}{rgb}{0.34,0.34,0.34}
\definecolor{gray35}{rgb}{0.35,0.35,0.35}
\definecolor{gray36}{rgb}{0.36,0.36,0.36}
\definecolor{gray37}{rgb}{0.37,0.37,0.37}
\definecolor{gray38}{rgb}{0.38,0.38,0.38}
\definecolor{gray39}{rgb}{0.39,0.39,0.39}
\definecolor{gray3}{rgb}{0.03,0.03,0.03}
\definecolor{gray40}{rgb}{0.40,0.40,0.40}
\definecolor{gray41}{rgb}{0.41,0.41,0.41}
\definecolor{gray42}{rgb}{0.42,0.42,0.42}
\definecolor{gray43}{rgb}{0.43,0.43,0.43}
\definecolor{gray44}{rgb}{0.44,0.44,0.44}
\definecolor{gray45}{rgb}{0.45,0.45,0.45}
\definecolor{gray46}{rgb}{0.46,0.46,0.46}
\definecolor{gray47}{rgb}{0.47,0.47,0.47}
\definecolor{gray48}{rgb}{0.48,0.48,0.48}
\definecolor{gray49}{rgb}{0.49,0.49,0.49}
\definecolor{gray4}{rgb}{0.04,0.04,0.04}
\definecolor{gray50}{rgb}{0.50,0.50,0.50}
\definecolor{gray51}{rgb}{0.51,0.51,0.51}
\definecolor{gray52}{rgb}{0.52,0.52,0.52}
\definecolor{gray53}{rgb}{0.53,0.53,0.53}
\definecolor{gray54}{rgb}{0.54,0.54,0.54}
\definecolor{gray55}{rgb}{0.55,0.55,0.55}
\definecolor{gray56}{rgb}{0.56,0.56,0.56}
\definecolor{gray57}{rgb}{0.57,0.57,0.57}
\definecolor{gray58}{rgb}{0.58,0.58,0.58}
\definecolor{gray59}{rgb}{0.59,0.59,0.59}
\definecolor{gray5}{rgb}{0.05,0.05,0.05}
\definecolor{gray60}{rgb}{0.60,0.60,0.60}
\definecolor{gray61}{rgb}{0.61,0.61,0.61}
\definecolor{gray62}{rgb}{0.62,0.62,0.62}
\definecolor{gray63}{rgb}{0.63,0.63,0.63}
\definecolor{gray64}{rgb}{0.64,0.64,0.64}
\definecolor{gray65}{rgb}{0.65,0.65,0.65}
\definecolor{gray66}{rgb}{0.66,0.66,0.66}
\definecolor{gray67}{rgb}{0.67,0.67,0.67}
\definecolor{gray68}{rgb}{0.68,0.68,0.68}
\definecolor{gray69}{rgb}{0.69,0.69,0.69}
\definecolor{gray6}{rgb}{0.06,0.06,0.06}
\definecolor{gray70}{rgb}{0.70,0.70,0.70}
\definecolor{gray71}{rgb}{0.71,0.71,0.71}
\definecolor{gray72}{rgb}{0.72,0.72,0.72}
\definecolor{gray73}{rgb}{0.73,0.73,0.73}
\definecolor{gray74}{rgb}{0.74,0.74,0.74}
\definecolor{gray75}{rgb}{0.75,0.75,0.75}
\definecolor{gray76}{rgb}{0.76,0.76,0.76}
\definecolor{gray77}{rgb}{0.77,0.77,0.77}
\definecolor{gray78}{rgb}{0.78,0.78,0.78}
\definecolor{gray79}{rgb}{0.79,0.79,0.79}
\definecolor{gray7}{rgb}{0.07,0.07,0.07}
\definecolor{gray80}{rgb}{0.80,0.80,0.80}
\definecolor{gray81}{rgb}{0.81,0.81,0.81}
\definecolor{gray82}{rgb}{0.82,0.82,0.82}
\definecolor{gray83}{rgb}{0.83,0.83,0.83}
\definecolor{gray84}{rgb}{0.84,0.84,0.84}
\definecolor{gray85}{rgb}{0.85,0.85,0.85}
\definecolor{gray86}{rgb}{0.86,0.86,0.86}
\definecolor{gray87}{rgb}{0.87,0.87,0.87}
\definecolor{gray88}{rgb}{0.88,0.88,0.88}
\definecolor{gray89}{rgb}{0.89,0.89,0.89}
\definecolor{gray8}{rgb}{0.08,0.08,0.08}
\definecolor{gray90}{rgb}{0.90,0.90,0.90}
\definecolor{gray91}{rgb}{0.91,0.91,0.91}
\definecolor{gray92}{rgb}{0.92,0.92,0.92}
\definecolor{gray93}{rgb}{0.93,0.93,0.93}
\definecolor{gray94}{rgb}{0.94,0.94,0.94}
\definecolor{gray95}{rgb}{0.95,0.95,0.95}
\definecolor{gray96}{rgb}{0.96,0.96,0.96}
\definecolor{gray97}{rgb}{0.97,0.97,0.97}
\definecolor{gray98}{rgb}{0.98,0.98,0.98}
\definecolor{gray99}{rgb}{0.99,0.99,0.99}
\definecolor{gray9}{rgb}{0.09,0.09,0.09}
\definecolor{gray}{rgb}{0.75,0.75,0.75}
\definecolor{green1}{rgb}{0.00,1.00,0.00}
\definecolor{green2}{rgb}{0.00,0.93,0.00}
\definecolor{green3}{rgb}{0.00,0.80,0.00}
\definecolor{green4}{rgb}{0.00,0.55,0.00}
\definecolor{greenyellow}{rgb}{0.68,1.00,0.18}
\definecolor{green}{rgb}{0.00,1.00,0.00}
\definecolor{grey0}{rgb}{0.00,0.00,0.00}
\definecolor{grey100}{rgb}{1.00,1.00,1.00}
\definecolor{grey10}{rgb}{0.10,0.10,0.10}
\definecolor{grey11}{rgb}{0.11,0.11,0.11}
\definecolor{grey12}{rgb}{0.12,0.12,0.12}
\definecolor{grey13}{rgb}{0.13,0.13,0.13}
\definecolor{grey14}{rgb}{0.14,0.14,0.14}
\definecolor{grey15}{rgb}{0.15,0.15,0.15}
\definecolor{grey16}{rgb}{0.16,0.16,0.16}
\definecolor{grey17}{rgb}{0.17,0.17,0.17}
\definecolor{grey18}{rgb}{0.18,0.18,0.18}
\definecolor{grey19}{rgb}{0.19,0.19,0.19}
\definecolor{grey1}{rgb}{0.01,0.01,0.01}
\definecolor{grey20}{rgb}{0.20,0.20,0.20}
\definecolor{grey21}{rgb}{0.21,0.21,0.21}
\definecolor{grey22}{rgb}{0.22,0.22,0.22}
\definecolor{grey23}{rgb}{0.23,0.23,0.23}
\definecolor{grey24}{rgb}{0.24,0.24,0.24}
\definecolor{grey25}{rgb}{0.25,0.25,0.25}
\definecolor{grey26}{rgb}{0.26,0.26,0.26}
\definecolor{grey27}{rgb}{0.27,0.27,0.27}
\definecolor{grey28}{rgb}{0.28,0.28,0.28}
\definecolor{grey29}{rgb}{0.29,0.29,0.29}
\definecolor{grey2}{rgb}{0.02,0.02,0.02}
\definecolor{grey30}{rgb}{0.30,0.30,0.30}
\definecolor{grey31}{rgb}{0.31,0.31,0.31}
\definecolor{grey32}{rgb}{0.32,0.32,0.32}
\definecolor{grey33}{rgb}{0.33,0.33,0.33}
\definecolor{grey34}{rgb}{0.34,0.34,0.34}
\definecolor{grey35}{rgb}{0.35,0.35,0.35}
\definecolor{grey36}{rgb}{0.36,0.36,0.36}
\definecolor{grey37}{rgb}{0.37,0.37,0.37}
\definecolor{grey38}{rgb}{0.38,0.38,0.38}
\definecolor{grey39}{rgb}{0.39,0.39,0.39}
\definecolor{grey3}{rgb}{0.03,0.03,0.03}
\definecolor{grey40}{rgb}{0.40,0.40,0.40}
\definecolor{grey41}{rgb}{0.41,0.41,0.41}
\definecolor{grey42}{rgb}{0.42,0.42,0.42}
\definecolor{grey43}{rgb}{0.43,0.43,0.43}
\definecolor{grey44}{rgb}{0.44,0.44,0.44}
\definecolor{grey45}{rgb}{0.45,0.45,0.45}
\definecolor{grey46}{rgb}{0.46,0.46,0.46}
\definecolor{grey47}{rgb}{0.47,0.47,0.47}
\definecolor{grey48}{rgb}{0.48,0.48,0.48}
\definecolor{grey49}{rgb}{0.49,0.49,0.49}
\definecolor{grey4}{rgb}{0.04,0.04,0.04}
\definecolor{grey50}{rgb}{0.50,0.50,0.50}
\definecolor{grey51}{rgb}{0.51,0.51,0.51}
\definecolor{grey52}{rgb}{0.52,0.52,0.52}
\definecolor{grey53}{rgb}{0.53,0.53,0.53}
\definecolor{grey54}{rgb}{0.54,0.54,0.54}
\definecolor{grey55}{rgb}{0.55,0.55,0.55}
\definecolor{grey56}{rgb}{0.56,0.56,0.56}
\definecolor{grey57}{rgb}{0.57,0.57,0.57}
\definecolor{grey58}{rgb}{0.58,0.58,0.58}
\definecolor{grey59}{rgb}{0.59,0.59,0.59}
\definecolor{grey5}{rgb}{0.05,0.05,0.05}
\definecolor{grey60}{rgb}{0.60,0.60,0.60}
\definecolor{grey61}{rgb}{0.61,0.61,0.61}
\definecolor{grey62}{rgb}{0.62,0.62,0.62}
\definecolor{grey63}{rgb}{0.63,0.63,0.63}
\definecolor{grey64}{rgb}{0.64,0.64,0.64}
\definecolor{grey65}{rgb}{0.65,0.65,0.65}
\definecolor{grey66}{rgb}{0.66,0.66,0.66}
\definecolor{grey67}{rgb}{0.67,0.67,0.67}
\definecolor{grey68}{rgb}{0.68,0.68,0.68}
\definecolor{grey69}{rgb}{0.69,0.69,0.69}
\definecolor{grey6}{rgb}{0.06,0.06,0.06}
\definecolor{grey70}{rgb}{0.70,0.70,0.70}
\definecolor{grey71}{rgb}{0.71,0.71,0.71}
\definecolor{grey72}{rgb}{0.72,0.72,0.72}
\definecolor{grey73}{rgb}{0.73,0.73,0.73}
\definecolor{grey74}{rgb}{0.74,0.74,0.74}
\definecolor{grey75}{rgb}{0.75,0.75,0.75}
\definecolor{grey76}{rgb}{0.76,0.76,0.76}
\definecolor{grey77}{rgb}{0.77,0.77,0.77}
\definecolor{grey78}{rgb}{0.78,0.78,0.78}
\definecolor{grey79}{rgb}{0.79,0.79,0.79}
\definecolor{grey7}{rgb}{0.07,0.07,0.07}
\definecolor{grey80}{rgb}{0.80,0.80,0.80}
\definecolor{grey81}{rgb}{0.81,0.81,0.81}
\definecolor{grey82}{rgb}{0.82,0.82,0.82}
\definecolor{grey83}{rgb}{0.83,0.83,0.83}
\definecolor{grey84}{rgb}{0.84,0.84,0.84}
\definecolor{grey85}{rgb}{0.85,0.85,0.85}
\definecolor{grey86}{rgb}{0.86,0.86,0.86}
\definecolor{grey87}{rgb}{0.87,0.87,0.87}
\definecolor{grey88}{rgb}{0.88,0.88,0.88}
\definecolor{grey89}{rgb}{0.89,0.89,0.89}
\definecolor{grey8}{rgb}{0.08,0.08,0.08}
\definecolor{grey90}{rgb}{0.90,0.90,0.90}
\definecolor{grey91}{rgb}{0.91,0.91,0.91}
\definecolor{grey92}{rgb}{0.92,0.92,0.92}
\definecolor{grey93}{rgb}{0.93,0.93,0.93}
\definecolor{grey94}{rgb}{0.94,0.94,0.94}
\definecolor{grey95}{rgb}{0.95,0.95,0.95}
\definecolor{grey96}{rgb}{0.96,0.96,0.96}
\definecolor{grey97}{rgb}{0.97,0.97,0.97}
\definecolor{grey98}{rgb}{0.98,0.98,0.98}
\definecolor{grey99}{rgb}{0.99,0.99,0.99}
\definecolor{grey9}{rgb}{0.09,0.09,0.09}
\definecolor{grey}{rgb}{0.75,0.75,0.75}
\definecolor{honeydew1}{rgb}{0.94,1.00,0.94}
\definecolor{honeydew2}{rgb}{0.88,0.93,0.88}
\definecolor{honeydew3}{rgb}{0.76,0.80,0.76}
\definecolor{honeydew4}{rgb}{0.51,0.55,0.51}
\definecolor{honeydew}{rgb}{0.94,1.00,0.94}
\definecolor{hotpink}{rgb}{1.00,0.41,0.71}
\definecolor{indianred}{rgb}{0.80,0.36,0.36}
\definecolor{ivory1}{rgb}{1.00,1.00,0.94}
\definecolor{ivory2}{rgb}{0.93,0.93,0.88}
\definecolor{ivory3}{rgb}{0.80,0.80,0.76}
\definecolor{ivory4}{rgb}{0.55,0.55,0.51}
\definecolor{ivory}{rgb}{1.00,1.00,0.94}
\definecolor{khaki1}{rgb}{1.00,0.96,0.56}
\definecolor{khaki2}{rgb}{0.93,0.90,0.52}
\definecolor{khaki3}{rgb}{0.80,0.78,0.45}
\definecolor{khaki4}{rgb}{0.55,0.53,0.31}
\definecolor{khaki}{rgb}{0.94,0.90,0.55}
\definecolor{lavenderblush}{rgb}{1.00,0.94,0.96}
\definecolor{lavender}{rgb}{0.90,0.90,0.98}
\definecolor{lawngreen}{rgb}{0.49,0.99,0.00}
\definecolor{lemonchiffon}{rgb}{1.00,0.98,0.80}
\definecolor{lightblue}{rgb}{0.68,0.85,0.90}
\definecolor{lightcoral}{rgb}{0.94,0.50,0.50}
\definecolor{lightcyan}{rgb}{0.88,1.00,1.00}
\definecolor{lightgoldenrod}{rgb}{0.93,0.87,0.51}
\definecolor{lightgoldenrod}{rgb}{0.98,0.98,0.82}
\definecolor{lightgray}{rgb}{0.83,0.83,0.83}
\definecolor{lightgreen}{rgb}{0.56,0.93,0.56}
\definecolor{lightgrey}{rgb}{0.83,0.83,0.83}
\definecolor{lightpink}{rgb}{1.00,0.71,0.76}
\definecolor{lightsalmon}{rgb}{1.00,0.63,0.48}
\definecolor{lightsea}{rgb}{0.13,0.70,0.67}
\definecolor{lightsky}{rgb}{0.53,0.81,0.98}
\definecolor{lightslate}{rgb}{0.47,0.53,0.60}
\definecolor{lightslate}{rgb}{0.47,0.53,0.60}
\definecolor{lightslate}{rgb}{0.52,0.44,1.00}
\definecolor{lightsteel}{rgb}{0.69,0.77,0.87}
\definecolor{lightyellow}{rgb}{1.00,1.00,0.88}
\definecolor{limegreen}{rgb}{0.20,0.80,0.20}
\definecolor{linen}{rgb}{0.98,0.94,0.90}
\definecolor{magenta1}{rgb}{1.00,0.00,1.00}
\definecolor{magenta2}{rgb}{0.93,0.00,0.93}
\definecolor{magenta3}{rgb}{0.80,0.00,0.80}
\definecolor{magenta4}{rgb}{0.55,0.00,0.55}
\definecolor{magenta}{rgb}{1.00,0.00,1.00}
\definecolor{maroon1}{rgb}{1.00,0.20,0.70}
\definecolor{maroon2}{rgb}{0.93,0.19,0.65}
\definecolor{maroon3}{rgb}{0.80,0.16,0.56}
\definecolor{maroon4}{rgb}{0.55,0.11,0.38}
\definecolor{maroon}{rgb}{0.69,0.19,0.38}
\definecolor{mediumaquamarine}{rgb}{0.40,0.80,0.67}
\definecolor{mediumblue}{rgb}{0.00,0.00,0.80}
\definecolor{mediumorchid}{rgb}{0.73,0.33,0.83}
\definecolor{mediumpurple}{rgb}{0.58,0.44,0.86}
\definecolor{mediumsea}{rgb}{0.24,0.70,0.44}
\definecolor{mediumslate}{rgb}{0.48,0.41,0.93}
\definecolor{mediumspring}{rgb}{0.00,0.98,0.60}
\definecolor{mediumturquoise}{rgb}{0.28,0.82,0.80}
\definecolor{mediumviolet}{rgb}{0.78,0.08,0.52}
\definecolor{midnightblue}{rgb}{0.10,0.10,0.44}
\definecolor{mintcream}{rgb}{0.96,1.00,0.98}
\definecolor{mistyrose}{rgb}{1.00,0.89,0.88}
\definecolor{moccasin}{rgb}{1.00,0.89,0.71}
\definecolor{navajowhite}{rgb}{1.00,0.87,0.68}
\definecolor{navyblue}{rgb}{0.00,0.00,0.50}
\definecolor{navy}{rgb}{0.00,0.00,0.50}
\definecolor{oldlace}{rgb}{0.99,0.96,0.90}
\definecolor{olivedrab}{rgb}{0.42,0.56,0.14}
\definecolor{orange1}{rgb}{1.00,0.65,0.00}
\definecolor{orange2}{rgb}{0.93,0.60,0.00}
\definecolor{orange3}{rgb}{0.80,0.52,0.00}
\definecolor{orange4}{rgb}{0.55,0.35,0.00}
\definecolor{orangered}{rgb}{1.00,0.27,0.00}
\definecolor{orange}{rgb}{1.00,0.65,0.00}
\definecolor{orchid1}{rgb}{1.00,0.51,0.98}
\definecolor{orchid2}{rgb}{0.93,0.48,0.91}
\definecolor{orchid3}{rgb}{0.80,0.41,0.79}
\definecolor{orchid4}{rgb}{0.55,0.28,0.54}
\definecolor{orchid}{rgb}{0.85,0.44,0.84}
\definecolor{palegoldenrod}{rgb}{0.93,0.91,0.67}
\definecolor{palegreen}{rgb}{0.60,0.98,0.60}
\definecolor{paleturquoise}{rgb}{0.69,0.93,0.93}
\definecolor{paleviolet}{rgb}{0.86,0.44,0.58}
\definecolor{papayawhip}{rgb}{1.00,0.94,0.84}
\definecolor{peachpuff}{rgb}{1.00,0.85,0.73}
\definecolor{peru}{rgb}{0.80,0.52,0.25}
\definecolor{pink1}{rgb}{1.00,0.71,0.77}
\definecolor{pink2}{rgb}{0.93,0.66,0.72}
\definecolor{pink3}{rgb}{0.80,0.57,0.62}
\definecolor{pink4}{rgb}{0.55,0.39,0.42}
\definecolor{pink}{rgb}{1.00,0.75,0.80}
\definecolor{plum1}{rgb}{1.00,0.73,1.00}
\definecolor{plum2}{rgb}{0.93,0.68,0.93}
\definecolor{plum3}{rgb}{0.80,0.59,0.80}
\definecolor{plum4}{rgb}{0.55,0.40,0.55}
\definecolor{plum}{rgb}{0.87,0.63,0.87}
\definecolor{powderblue}{rgb}{0.69,0.88,0.90}
\definecolor{purple1}{rgb}{0.61,0.19,1.00}
\definecolor{purple2}{rgb}{0.57,0.17,0.93}
\definecolor{purple3}{rgb}{0.49,0.15,0.80}
\definecolor{purple4}{rgb}{0.33,0.10,0.55}
\definecolor{purple}{rgb}{0.63,0.13,0.94}
\definecolor{red1}{rgb}{1.00,0.00,0.00}
\definecolor{red2}{rgb}{0.93,0.00,0.00}
\definecolor{red3}{rgb}{0.80,0.00,0.00}
\definecolor{red4}{rgb}{0.55,0.00,0.00}
\definecolor{red}{rgb}{1.00,0.00,0.00}
\definecolor{rosybrown}{rgb}{0.74,0.56,0.56}
\definecolor{royalblue}{rgb}{0.25,0.41,0.88}
\definecolor{saddlebrown}{rgb}{0.55,0.27,0.07}
\definecolor{salmon1}{rgb}{1.00,0.55,0.41}
\definecolor{salmon2}{rgb}{0.93,0.51,0.38}
\definecolor{salmon3}{rgb}{0.80,0.44,0.33}
\definecolor{salmon4}{rgb}{0.55,0.30,0.22}
\definecolor{salmon}{rgb}{0.98,0.50,0.45}
\definecolor{sandybrown}{rgb}{0.96,0.64,0.38}
\definecolor{seagreen}{rgb}{0.18,0.55,0.34}
\definecolor{seashell1}{rgb}{1.00,0.96,0.93}
\definecolor{seashell2}{rgb}{0.93,0.90,0.87}
\definecolor{seashell3}{rgb}{0.80,0.77,0.75}
\definecolor{seashell4}{rgb}{0.55,0.53,0.51}
\definecolor{seashell}{rgb}{1.00,0.96,0.93}
\definecolor{sienna1}{rgb}{1.00,0.51,0.28}
\definecolor{sienna2}{rgb}{0.93,0.47,0.26}
\definecolor{sienna3}{rgb}{0.80,0.41,0.22}
\definecolor{sienna4}{rgb}{0.55,0.28,0.15}
\definecolor{sienna}{rgb}{0.63,0.32,0.18}
\definecolor{skyblue}{rgb}{0.53,0.81,0.92}
\definecolor{slateblue}{rgb}{0.42,0.35,0.80}
\definecolor{slategray}{rgb}{0.44,0.50,0.56}
\definecolor{slategrey}{rgb}{0.44,0.50,0.56}
\definecolor{snow1}{rgb}{1.00,0.98,0.98}
\definecolor{snow2}{rgb}{0.93,0.91,0.91}
\definecolor{snow3}{rgb}{0.80,0.79,0.79}
\definecolor{snow4}{rgb}{0.55,0.54,0.54}
\definecolor{snow}{rgb}{1.00,0.98,0.98}
\definecolor{springgreen}{rgb}{0.00,1.00,0.50}
\definecolor{steelblue}{rgb}{0.27,0.51,0.71}
\definecolor{tan1}{rgb}{1.00,0.65,0.31}
\definecolor{tan2}{rgb}{0.93,0.60,0.29}
\definecolor{tan3}{rgb}{0.80,0.52,0.25}
\definecolor{tan4}{rgb}{0.55,0.35,0.17}
\definecolor{tan}{rgb}{0.82,0.71,0.55}
\definecolor{thistle1}{rgb}{1.00,0.88,1.00}
\definecolor{thistle2}{rgb}{0.93,0.82,0.93}
\definecolor{thistle3}{rgb}{0.80,0.71,0.80}
\definecolor{thistle4}{rgb}{0.55,0.48,0.55}
\definecolor{thistle}{rgb}{0.85,0.75,0.85}
\definecolor{tomato1}{rgb}{1.00,0.39,0.28}
\definecolor{tomato2}{rgb}{0.93,0.36,0.26}
\definecolor{tomato3}{rgb}{0.80,0.31,0.22}
\definecolor{tomato4}{rgb}{0.55,0.21,0.15}
\definecolor{tomato}{rgb}{1.00,0.39,0.28}
\definecolor{turquoise1}{rgb}{0.00,0.96,1.00}
\definecolor{turquoise2}{rgb}{0.00,0.90,0.93}
\definecolor{turquoise3}{rgb}{0.00,0.77,0.80}
\definecolor{turquoise4}{rgb}{0.00,0.53,0.55}
\definecolor{turquoise}{rgb}{0.25,0.88,0.82}
\definecolor{violetred}{rgb}{0.82,0.13,0.56}
\definecolor{violet}{rgb}{0.93,0.51,0.93}
\definecolor{wheat1}{rgb}{1.00,0.91,0.73}
\definecolor{wheat2}{rgb}{0.93,0.85,0.68}
\definecolor{wheat3}{rgb}{0.80,0.73,0.59}
\definecolor{wheat4}{rgb}{0.55,0.49,0.40}
\definecolor{wheat}{rgb}{0.96,0.87,0.70}
\definecolor{whitesmoke}{rgb}{0.96,0.96,0.96}
\definecolor{white}{rgb}{1.00,1.00,1.00}
\definecolor{yellow1}{rgb}{1.00,1.00,0.00}
\definecolor{yellow2}{rgb}{0.93,0.93,0.00}
\definecolor{yellow3}{rgb}{0.80,0.80,0.00}
\definecolor{yellow4}{rgb}{0.55,0.55,0.00}
\definecolor{yellowgreen}{rgb}{0.60,0.80,0.20}
\definecolor{yellow}{rgb}{1.00,1.00,0.00}

\maketitle

\begin{abstract}

\end{abstract}

In this paper we prove a homogenous generalization of the lambda determinant formula of Mills, Robbins and Rumsey. In our formula the parameters depends on two indices. Our result also extends a recent formula of Di Francesco.

\section{Introduction}

An \emph{alternating sign matrix} is a square matrix of $0$'s $1$'s and $-1$'s such that the sum of each row and column is $1$ and the non-zero entries in each row and column alternate in sign. For example:
\[ X =
\left ( \begin{matrix}
0 & 1 & 0 & 0  \\
1 & -1 & 1 & 0  \\
0 & 1 & -1 & 1  \\
0 & 0 & 1 & 0 
\end{matrix} \right )
\]

Alternating sign matrices arise naturally in Dodgson's condensation method for calculating $\lambda$-determinants \cite{asm}.

For each $k=0..n$ let us denote by $x[k]$ the doubly indexed collection of variables $x[k]_{i,j}$ with indices running from $i,j = 1..(n-k+1)$. One should think of these variables as forming a square pyramid with base of dimension $n+1$ by $n+1$. The index $k$ determines the ``height'' 
of the variable in the pyramid.

The variables $x[0]$ and $x[1]$ are to be thought of as initial conditions. The remaining $x[k]$ are defined in terms of the following \emph{octahedral recurrence}:
\begin{align} \label{rec1}
x[k+1]_{i,j} =  \frac{ \mu_{i,n-k-j} x[k]_{i,j} x[k]_{i+1,j+1} + \lambda_{i,j} x[k]_{i,j+1} x[k]_{i+1,j}}{x[k-1]_{i+1,j+1}}
\end{align}

The main result of this paper is a closed form expression for $x[k]_{1,1}$.
Our result generalizes the result obtained by  Di Francesco \cite{cluster}, who considered
coefficients $\lambda_{ij}\equiv\lambda_{i-j}$ and $\mu_{ij}\equiv\mu_{i-j}$.


The outline of this paper is as follows. We begin with some definitions which are necessary in order to write down the closed form expression. In section \ref{up-down} we introduce \emph{left cumulant matrices} and a pair of \emph{up / down operators}. In section  \ref{more} we introduce \emph{right cumulant matrices} and a second pair of up / down operators which are closely related to the first. Finally in section \ref{proof} we prove our main theorem.

\section{closed form expression}\label{closed}

For each $n$ by $n$ alternating sign matrix $B$ let $\overline{B}$ be the matrix whose $(i,j)$-th entry is equal to the sum of the entries lying above and to the left of the $(i,j)$-th entry of $B$. Similarly, let $\underline{B}
$ be the matrix whose $(i,j)$-th entry is equal to the sum of the entries lying above and to the \emph{right} of the $(i,j)$-th entry of $B$. For example:
\[ X =
\left ( \begin{matrix}
0 & 1 & 0 & 0  \\
1 & -1 & 1 & 0  \\
0 & 1 & -1 & 1  \\
0 & 0 & 1 & 0 
\end{matrix} \right )
\quad
 \overline{X} =
\left ( \begin{matrix}
0 & 1 & 1 & 1  \\
1 & 1 & 2 & 2  \\
1 & 2 & 2 & 3  \\
1 & 2 & 3 & 4 
\end{matrix} \right )
\quad
 \underline{X} =
\left ( \begin{matrix}
1 & 1 & 0 & 0  \\
2 & 1 & 1 & 0  \\
3 & 2 & 1 & 1  \\
4 & 3 & 2 & 1 
\end{matrix} \right )
\]

We shall refer to $\overline{B}$ as the \emph{left cumulant} matrix of $B$ and $\underline{B}$ as the \emph{right cumulant} matrix of $B$.
The original alternating sign matrix may be recovered by the formula:

\begin{align} 
B_{ij} & = \overline{B}_{ij} + \overline{B}_{i-1,j-1} - \overline{B}_{i,j-1} - \overline{B}_{i-1,j} \label{seconddiff}\\
& = \underline{B}_{ij} + \underline{B}_{i-1,j+1} - \underline{B}_{i,j+1} - \underline{B}_{i-1,j}\label{seconddiff2}
\end{align}
If the indices are out of range, then the value of $B_{ij}$ is taken to be zero.

\begin{lemma}\label{reflect}
If $B'$ is the alternating sign matrix obtained from $B$ by multiplying on the right by the maximum permutation then 
$\underline{B'}$ is the matrix obtained from $\overline{B}$ by multiplying on the right by the maximum permutation.
\end{lemma}
We shall make use of the notation:
\[ \underline{B} = (\overline{B'})' \]

\begin{lemma}\label{lemma}
For all $i,j$ we have:
\[\overline{B}_{i,j} + \underline{B}_{i,j+1} = i \]
\begin{proof}
The left hand side is equal to the sum of all the entries of the alternating sign matrix $B$ in the first $i$ rows. Since the sum of entries in each row of $B$ is equal to $1$, the final result is equal to $i$ as claimed. 
\end{proof}
\end{lemma}

Comparing matrices entrywise,
the $\overline{B}$ of size $n$ form a lattice. We remark that this lattice coincides with the completion of the Bruhat order to alternating
sign matrices as carried out in Lascoux and
Sch\"utzenberger \cite{treillis}.
One can apply the same operation with \underline{B} to form a dual lattice.

Let us define the \emph{lambda weight} of a $k$ by $k$ alternating sign matrix $B$ to be:
\[ \lambda^{F(B)} = \lambda^{\overline{I} - \overline{B}} = \prod_{i,j=1}^k \lambda_{i,j}^{\min(i,j) - \overline{B}_{i,j}} \]
Similarly, let us define the \emph{mu weight} of an $k$ by $k$ alternating sign matrix $B$ to be:
\[ \mu^{G(B)} = \mu^{\underline{I} - \underline{B}} = \prod_{i,j=1}^k \mu_{i,j}^{\max(i-j+1,0) - \underline{B}_{i,k+1-j}} \]
The \emph{standard weight} of an alternating sign matrix $B$ is simply:
\[ x^B = \prod_{i,j=1}^k x_{i,j}^{B_{i,j}} \]

Robbins and Rumsey \cite{robrum} defined two multiplicity free operators acting on the vector space spanned by alternating sign matrices, which we shall discuss in section \ref{up-down}:

\begin{align*}
\mathfrak{U} : \ASM(n) \to \mathbb{Z} [\ASM(n+1)] \\
\mathfrak{D} : \ASM(n) \to \mathbb{Z} [\ASM(n-1)]
\end{align*}


Our closed form expression for $x[k]_{1,1}$ now takes the form;

\begin{equation} \label{double}
\boxed{x[k]_{1,1} = \sum_{\substack{(A,B) \\ |B| = k, |A| = k-1 \\ A \in \mathfrak{D}(B)}} 
\lambda^{F(B)} s(\lambda)^{-F(A)} \mu^{G(B)} t(\mu)^{-G(A)}x[1]^B s(x[0])^{-A} }
\end{equation}
where: 
\begin{align}
s(z)_{i,j} & = z_{i+1,j+1} \\
t(z)_{i,j} & = z_{i+1,j-1}
\end{align}

Note that this formula shows that $x[k]_{1,1}$ is a Laurent polynomial, and not just a rational function as would be expected from its recursive definition. This is an exanple of te Laurent phenomenon. See, for example \cite{fomin}.

\section{Up and Down operators}\label{up-down}

We shall now define the multiplicity free operators acting on the vector space spanned by alternating sign matrices mentioned in the previous section:

\begin{align*}
\mathfrak{U} : \ASM(n) \to \mathbb{Z} [\ASM(n+1)] \\
\mathfrak{D} : \ASM(n) \to \mathbb{Z} [\ASM(n-1)]
\end{align*}

These operators have the property that $B \in \ASM(n)$ contains $r$ ones and $s$ negative ones then number of terms occuring in  $\mathfrak{U}(B)$ is $2^r$ while the number of terns occuring in $\mathfrak{D}(B)$ is $2^s$.

If we fix an order on the $-1$'s of $B$ then each element $A$ of $\mathfrak{D}(B)$ is naturally indexed by a binary string.
Similarly if we fix an order of the $1$'s in $B$ then each element $C$ of $\mathfrak{U}(B)$  is indexed by a binary string.

To define these operators we shall need the notion of \emph{left interlacing matrices}: 

\[
\left ( \begin{matrix}
\overline{B}_{1,1} & & \overline{B}_{1,2} & & \overline{B}_{1,3}& & \overline{B}_{1,4}  \\
 &  {\color{red}\overline{A}_{1,1}}&  & \color{red}\overline{A}_{1,2} &  & \color{red}\overline{A}_{1,3} &   \\
\overline{B}_{2,1} & & \overline{B}_{2,2} & & \overline{B}_{2,3} & & \overline{B}_{2,4} \\
 &  {\color{red}\overline{A}_{2,1}}&  & \color{red}\overline{A}_{2,2} &  & \color{red}\overline{A}_{2,3} &   \\
\overline{B}_{3,1} & & \overline{B}_{3,2} & & \overline{B}_{3,3} & & \overline{B}_{3,4}  \\
 &  {\color{red}\overline{A}_{3,1}}&  & \color{red}\overline{A}_{3,2} &  & \color{red}\overline{A}_{3,3} &   \\
\overline{B}_{4,1} & & \overline{B}_{4,2} & & \overline{B}_{4,3} & & \overline{B}_{4,4} 
\end{matrix} \right )
\]

The conditions on the matrix $\overline{A}$ are as follows:

\[
\left ( \begin{matrix}
x & & y \\
 &  {\color{red}a}& \\
z & & w 
\end{matrix} \right )
\qquad \qquad
\boxed{x,w-1 \leq {\color{red}a} \leq y,z}
\]

An example:

\[
\left ( \begin{matrix}
0 & & 1 & & 1 & & 1  \\
 &  {\color{red}\{0,1\}}&  & \color{red}1 &  & \color{red}1 &   \\
1 & & 1 & & 2 & & 2  \\
 &  \color{red}1&  & \color{red}\{1,2\}&  & \color{red}2 &   \\
1 & & 2 & & 2 & & 3  \\
 & \color{red} 1&  & \color{red}2 &  & \color{red}3 &   \\
1 & & 2 & & 3 & & 4 
\end{matrix} \right )
\]

Above and to the left of a $-1$ in the alternating sign matrix $B$ there are two possible choices for the corresponding value of the left cumulant matrix $\overline{A}$. At all other positions there is a single choice \cite{robrum}.



\begin{center}
\begin{tikzpicture}
\begin{scope}
\path (5,0) node[](min){
$
\overline{A}_{00} = \left ( \begin{matrix}
 0 & 1 & 1 \\
1 & 1 & 2 \\
1 & 2 & 3 \\
\end{matrix} \right ) 
$
};

\path (2,2) node[](left){
$
\overline{A}_{01} = \left ( \begin{matrix}
 0 & 1 & 1 \\
1 & 2 & 2 \\
1 & 2 & 3 \\
\end{matrix} \right ) 
$
};
\path (8,2) node[](right){
$
\overline{A}_{10} = \left ( \begin{matrix}
 1 & 1 & 1 \\
1 & 1 & 2 \\
1 & 2 & 3 \\
\end{matrix} \right ) 
$
};

\path (5,4) node[](max){
$
\overline{A}_{11} = \left ( \begin{matrix}
 1 & 1 & 1 \\
1 & 2 & 2 \\
1 & 2 & 3 \\
\end{matrix} \right ) 
$
};

\draw (min) -- (left) -- (max) -- (right) -- (min);

\end{scope}
\end{tikzpicture}
\end{center}

Here are the corresponding alternating sign matrices:

\begin{center}
\begin{tikzpicture}
\begin{scope}
\path (5,0) node[](min){
$
A_{00} = \left ( \begin{matrix}
 0 & 1 & 0 \\
1 & -1 & 1 \\
0 & 1 & 0 \\
\end{matrix} \right )
$
};

\path (2,2) node[](left){
$
A_{01} = \left ( \begin{matrix}
 0 & 1 & 0 \\
1 & 0 & 0 \\
0 & 0 & 1 \\
\end{matrix} \right )
$
};
\path (8,2) node[](right){
$
A_{10} = \left ( \begin{matrix}
 1 & 0 & 0 \\
0 & 0 & 1 \\
0 & 1 & 0 \\
\end{matrix} \right )
$
};

\path (5,4) node[](max){
$
A_{11} = \left ( \begin{matrix}
 1 & 0 & 0 \\
0 & 1 & 0 \\
0 & 0 & 1 \\
\end{matrix} \right )
$
};

\draw (min) -- (left) -- (max) -- (right) -- (min);

\end{scope}
\end{tikzpicture}
\end{center}

One may check that adding one at position $(i,j)$ in the left cumulant matrix $\overline{A}$ is equivalent to adding the matrix 
$ \left (\begin{matrix} 1 & -1 \\ -1 & 1 \end{matrix} \right ) $  with upper left hand corner at position $(i,j)$
to the corresponding alternating sign matrix $A$ \cite{robrum}. 

In our example we have:
\[ \mathfrak{D}(X) = A_{0,0} + A_{1,0} + A_{0,1} + A_{1,1} \]

We shall be especially interested in the ``smallest'' matrix $A$ which is left interlacing with the matrix $B$ and which we denote by $A^{\min} = A_{00\cdots 0}$.
We have, by construction:
\[ \boxed{\overline{A}^{\min}_{ij} = \max(\overline{B}_{ij}, \overline{B}_{i+1,j+1} - 1)} \] 

The $\mathfrak{U}$ operator is defined similarly.

\[
\left ( \begin{matrix} 
\color{blue} \overline{C}_{1,1} & & \color{blue} \overline{C}_{1,2}  & & \color{blue} \overline{C}_{1,3}  \\
& \overline{B}_{1,1} & & \overline{B}_{1,2} & \\
\color{blue} \overline{C}_{2,1} & & \color{blue} \overline{C}_{2,2}  & & \color{blue} \overline{C}_{2,3}  \\
& \overline{B}_{2,1} & & \overline{B}_{2,2} & \\
\color{blue} \overline{C}_{3,1} & & \color{blue} \overline{C}_{3,2}  & & \color{blue} \overline{C}_{3,3}  \\
\end{matrix} \right )
\]

The rule for constructing all possible matrices $\overline{C}$ for a given matrix $\overline{B}$ is the last row and last column must be strictly increasing from $1$ to $n$
as well as that:

\[
\left ( \begin{matrix}
x & & y \\
 &  {\color{blue}c}& \\
z & & w 
\end{matrix} \right )
\qquad \qquad
\boxed{
y,z \leq {\color{blue}c} \leq w,x+1}
\]

Here is an example:
\[ Y = \left ( \begin{matrix} 0 & 1 \\ 1 & 0 \end{matrix} \right) \qquad 
\overline{Y} = \left ( \begin{matrix} 0 & 1 \\ 1 & 2 \end{matrix} \right ) \]
Interlacing matrices:
\[
\left ( \begin{matrix} 
\color{blue} 0 & & \color{blue} \{0,1\} & & \color{blue} 1 \\
& 0 & & 1 & \\
\color{blue} \{0,1\} & & \color{blue} 1 & & \color{blue} 2 \\
& 1 & & 2 & \\
\color{blue} 1 & & \color{blue} 2 & & \color{blue} 3 \\
\end{matrix} \right )
\]

Above and to the left of a $1$ in the alternating sign matrix $B$ there are two possible choices for the corresponding value of $\overline{C}$. At all other positions there is a single choice \cite{robrum}. 


\begin{center}
\begin{tikzpicture}
\begin{scope}
\path (5,4) node[](max){
$
\overline{C}_{11} = \left ( \begin{matrix}
 0 & 1 & 1 \\
1 & 1 & 2 \\
1 & 2 & 3 \\
\end{matrix} \right )
$
};

\path (2,2) node[](left){
$
\overline{C}_{01} = \left ( \begin{matrix}
 0 & 1 & 1 \\
0 & 1 & 2 \\\hfill \break
1 & 2 & 3 \\
\end{matrix} \right ) 
$
};
\path (8,2) node[](right){
$
\overline{C}_{10} = \left ( \begin{matrix}
 0 & 0 & 1 \\
1 & 1 & 2 \\
1 & 2 & 3 \\
\end{matrix} \right ) 
$
};

\path (5,0) node[](min){
$
\overline{C}_{00} = \left ( \begin{matrix}
 0 & 0 & 1 \\
0 & 1 & 2 \\
1 & 2 & 3 \\
\end{matrix} \right ) 
$
};

\draw (min) -- (left) -- (max) -- (right) -- (min);

\end{scope}
\end{tikzpicture}
\end{center}

Here are the corresponding alternating sign matrices:

\begin{center}
\begin{tikzpicture}
\begin{scope}
\path (5,0) node[](min){
$
C_{00} = \left ( \begin{matrix}
 0 & 0 & 1 \\
0 & 1 & 0 \\
1 & 0 & 0 \\
\end{matrix} \right ) 
$
};

\path (2,2) node[](left){
$
C_{01} = \left ( \begin{matrix}
 0 & 1 & 0 \\
0 & 0 & 1 \\
1 & 0 & 0 \\
\end{matrix} \right )
$
};
\path (8,2) node[](right){
$
C_{10} = \left ( \begin{matrix}
 0 & 0 & 1 \\
1 & 0 & 0 \\
0 & 1 & 0 \\
\end{matrix} \right )
$
};

\path (5,4) node[](max){
$
C_{11} = \left ( \begin{matrix}
 0 & 1 & 0 \\
1 & -1 & 1 \\
0 & 1 & 0 \\
\end{matrix} \right )
$
};

\draw (min) -- (left) -- (max) -- (right) -- (min);

\end{scope}
\end{tikzpicture}
\end{center}

One may check that, as expected, subtracting one at position $(i,j)$ in the left cumulant matrix $\overline{C}$ is equivalent to subtracting the matrix 
$ \left (\begin{matrix} 1 & -1 \\ -1 & 1 \end{matrix} \right ) $
 with upper left hand corner at position $(i,j)$ from the corresponding alternating sign matrix $C$ \cite{robrum}. 

In our example we have:
\[ \mathfrak{U}(Y) = C_{0,0} + C_{0,1} + C_{1,0} + C_{1,1} \]

This time we shall be especially interested in the ``largest'' matrix $C$ which is left interlacing with $B$ and which we denote by $C^{\max} = C_{11 \cdots 1}$.
We have, by construction:
\[ \boxed{\overline{C}^{\max}_{ij} = \min(\overline{B}_{ij}, \overline{B}_{i-1,j-1} + 1)} \] 

\section{More up-down operators}\label{more}

We will need a second set of up and down operators which are closely related to the first.

\begin{align*}
\mathfrak{U}^* : \ASM(n) \to \mathbb{Z} [\ASM(n+1)] \\
\mathfrak{D}^* : \ASM(n) \to \mathbb{Z} [\ASM(n-1)]
\end{align*}

To define these operators we make use of \emph{right interlacing matrices}:  

\[
\left ( \begin{matrix}
\underline{B}_{1,1} & & \underline{B}_{1,2} & & \underline{B}_{1,3}& & \underline{B}_{1,4}  \\
 &  {\color{red}\underline{A}^*_{1,1}}&  & \color{red}\underline{A}^*_{1,2} &  & \color{red}\underline{A}^*_{1,3} &   \\
\underline{B}_{2,1} & & \underline{B}_{2,2} & & \underline{B}_{2,3} & & \underline{B}_{2,4} \\
 &  {\color{red}\underline{A}^*_{2,1}}&  & \color{red}\underline{A}^*_{2,2} &  & \color{red}\underline{A}^*_{2,3} &   \\
\underline{B}_{3,1} & & \underline{B}_{3,2} & & \underline{B}_{3,3} & & \underline{B}_{3,4}  \\
 &  {\color{red}\underline{A}^*_{3,1}}&  & \color{red}\underline{A}^*_{3,2} &  & \color{red}\underline{A}^*_{3,3} &   \\
\underline{B}_{4,1} & & \underline{B}_{4,2} & & \underline{B}_{4,3} & & \underline{B}_{4,4} 
\end{matrix} \right )
\]

In the right interlacing case, the conditions on the matrix $\underline{A}^*$ are:

\[
\left ( \begin{matrix}
x & & y \\
 &  {\color{red}a}& \\
z & & w 
\end{matrix} \right )
\qquad \qquad
\boxed{
y,z-1 \leq {\color{red}a} \leq x,w}
\]

Continuing with our example matrix $X$:

\[
\left ( \begin{matrix}
1 & & 1 & & 0 & & 0  \\
 &  {\color{red}1}&  & \color{red}\{0,1\} &  & \color{red}0 &   \\
1 & & 1 & & 2 & & 2  \\
 &  \color{red}2&  & \color{red}1&  & \color{red}\{0,1\} &   \\
2 & & 1 & & 1 & & 0  \\
 & \color{red} 3&  & \color{red}2 &  & \color{red}1 &   \\
4 & & 3 & & 2 & & 1 
\end{matrix} \right )
\]

Above and to the \emph{right} of a $-1$ in the alternating sign matrix $B$ there are two possible choices for the corresponding value of the right cumulant matrix $\underline{A}^*$. Again, if we fix an order on the $-1$'s of $B$ then each element $A^*$ of $\mathfrak{D}(B)$ is naturally indexed by a binary string determining the position in the right cumulant matrix $\underline{A}^*$ where the larger of the two possible values was chosen. 


\begin{center}
\begin{tikzpicture}
\begin{scope}
\path (5,0) node[](min){
$
\underline{A}^*_{00} = \left ( \begin{matrix}
 1 & 0 & 0 \\
2 & 1 & 1 \\
3 & 2 & 1 \\
\end{matrix} \right ) 

$
};

\path (2,2) node[](left){
$
\underline{A}^*_{01} = \left ( \begin{matrix}
 1 & 0 & 0 \\
2 & 1 & 1 \\
3 & 2 & 1 \\
\end{matrix} \right ) 
$
};
\path (8,2) node[](right){
$
\underline{A}^*_{10} = \left ( \begin{matrix}
 1 & 1 & 0 \\
2 & 1 & 0 \\
3 & 2 & 1 \\
\end{matrix} \right ) 
$
};

\path (5,4) node[](max){
$
\underline{A}^*_{11} = \left ( \begin{matrix}
 1 & 1 & 0\\
2 & 1 & 1 \\
3 & 2 & 1 \\
\end{matrix} \right ) 
$
};

\draw (min) -- (left) -- (max) -- (right) -- (min);

\end{scope}
\end{tikzpicture}
\end{center}

Here are the corresponding alternating sign matrices:

\begin{center}
\begin{tikzpicture}
\begin{scope}
\path (5,0) node[](min){
$
A^*_{00} = \left ( \begin{matrix}
 1 & 0 & 0 \\
0 & 1 & 0 \\
0 & 0 & 1 \\
\end{matrix} \right )

$
};

\path (2,2) node[](left){
$
A^*_{01} = \left ( \begin{matrix}
 1 & 0 & 0 \\
0 & 0 & 1 \\
0 & 1 & 0 \\
\end{matrix} \right )

$
};
\path (8,2) node[](right){
$
A^*_{10} = \left ( \begin{matrix}
 0 & 1 & 0 \\
1 & 0 & 0 \\
0 & 0 & 1 \\
\end{matrix} \right ) 
$
};

\path (5,4) node[](max){
$
A^*_{11} = \left ( \begin{matrix}
 0 & 1 & 0 \\
1 & -1 & 1 \\
0 & 1 & 0 \\
\end{matrix} \right ) 
$
};

\draw (min) -- (left) -- (max) -- (right) -- (min);

\end{scope}
\end{tikzpicture}
\end{center}

Adding one at position $(i,j)$ in the right cumulant matrix $\underline{A}^*$ is equivalent to adding the matrix 
$ \left (\begin{matrix} -1 & 1 \\ 1 & -1 \end{matrix} \right ) $ with upper \emph{right} hand corner at position $(i,j)$
to the alternating sign matrix $A^*$.

If $B$ is an $n$ by $n$ alternating sign matrix then $\mathfrak{D}^*(B)$ is the sum of all $n-1$ by $n-1$ alternating sign matrices $A^*$ such that $\underline{A}^*$ is right interlacing with $\underline{B}^*$.

\hfill \break

Now for the $\mathfrak{U}^*$ operator.

\[
\left ( \begin{matrix} 
\color{blue} \underline{C}^*_{1,1} & & \color{blue} \underline{C}^*_{1,2}  & & \color{blue} \underline{C}^*_{1,3}  \\
& \underline{B}_{1,1} & & \underline{B}_{1,2} & \\
\color{blue} \underline{C}^*_{2,1} & & \color{blue} \underline{C}^*_{2,2}  & & \color{blue} \underline{C}^*_{2,3}  \\
& \underline{B}_{2,1} & & \underline{B}_{2,2} & \\
\color{blue} \underline{C}^*_{3,1} & & \color{blue} \underline{C}^*_{3,2}  & & \color{blue} \underline{C}^*_{3,3}  \\
\end{matrix} \right )
\]

The rule for constructing all possible matrices $\underline{C}^*$ for a given matrix $\underline{B}$ is the last column must be strictly increasing from $1$ to $n+1$, the last row must be strictly decreasing from $n$ to $1$,
and:

\[
\left ( \begin{matrix}
x & & y \\
 &  {\color{blue}c}& \\
z & & w 
\end{matrix} \right )
\qquad \qquad
\boxed{
w,x \leq {\color{blue}c} \leq y+1,z}
\]

Here is an example:
\[ Y = \left ( \begin{matrix} 0 & 1 \\ 1 & 0 \end{matrix} \right) \qquad 
\underline{Y} = \left ( \begin{matrix} 1 & 1 \\ 2 & 1 \end{matrix} \right ) \]
Interlacing matrices:
\[
\left ( \begin{matrix} 
\color{blue} 1 & & \color{blue} 1 & & \color{blue} \{0,1\} \\
& 1 & & 1 & \\
\color{blue} 2 & & \color{blue} \{1,2\} & & \color{blue} 1 \\
& 2 & & 1 & \\
\color{blue} 3 & & \color{blue} 2 & & \color{blue} 1 \\
\end{matrix} \right )
\]

Above and to the \emph{right} of a $1$ in the alternating sign matrix $B$ there are two possible choices for the corresponding value of $\underline{C}^*$. At all other positions there is a single choice. 

Fixing an order on the $1$'s of $B$, for each element $C^*$ of $\mathfrak{U}^*(B)$ is naturally indexed by a binary string determing the position in the right cumulant matrix $\underline{C}^*$ where the larger of the two possible values was chosen.

\begin{center}
\begin{tikzpicture}
\begin{scope}
\path (5,4) node[](max){
$
\underline{C}^*_{11} = \left ( \begin{matrix}
 1 & 1 & 1 \\
2 & 2 & 1 \\
3 & 2 & 1 \\
\end{matrix} \right )
$
};

\path (2,2) node[](left){
$
\underline{C}^*_{01} = \left ( \begin{matrix}
 1 & 1 & 1 \\
2 & 1 & 1 \\\hfill \break
3 & 2 & 1 \\
\end{matrix} \right ) 
$
};
\path (8,2) node[](right){
$
\underline{C}^*_{10} = \left ( \begin{matrix}
 1 & 1 & 0 \\
2 & 2 & 1 \\
3 & 2 & 1 \\
\end{matrix} \right ) 
$
};

\path (5,0) node[](min){
$
\underline{C}^*_{00} = \left ( \begin{matrix}
 1 & 1 & 0 \\
2 & 1 & 1 \\
3 & 2 & 1 \\
\end{matrix} \right ) 
$
};

\draw (min) -- (left) -- (max) -- (right) -- (min);

\end{scope}
\end{tikzpicture}
\end{center}

Here are the corresponding alternating sign matrices:

\begin{center}
\begin{tikzpicture}
\begin{scope}
\path (5,0) node[](min){
$
C^*_{00} = \left ( \begin{matrix}
 0 & 1 & 0 \\
1 & -1 & 1 \\
0 & 1 & 0 \\
\end{matrix} \right )
$
};

\path (2,2) node[](left){
$
C^*_{01} = \left ( \begin{matrix}
 0 & 0 & 1 \\
1 & 0 & 0 \\
0 & 1 & 0 \\
\end{matrix} \right ) 
$
};
\path (8,2) node[](right){
$
C^*_{10} = \left ( \begin{matrix}
 0 & 1 & 0 \\
0 & 0 & 1 \\
1 & 0 & 0 \\
\end{matrix} \right )

$
};

\path (5,4) node[](max){
$
C^*_{11} = \left ( \begin{matrix}
 0 & 0 & 1 \\
0 & 1 & 0 \\
1 & 0 & 0 \\
\end{matrix} \right )  
$
};

\draw (min) -- (left) -- (max) -- (right) -- (min);

\end{scope}
\end{tikzpicture}
\end{center}

Subtracting one at position $(i,j)$ in the right cumulant matrix $\underline{C}^*$ is equivalent to subtracting the matrix 
$ \left (\begin{matrix} -1 & 1 \\ 1 & -1 \end{matrix} \right ) $
 with upper \emph{right} hand corner at position $(i,j)$ from the alternating sign matrix $C^*$.

\begin{proposition}
\[{A}^*_{\min} = A_{\max} \]
\begin{proof}
Consider the following segments of interlacing matrices:
\[
\left ( \begin{matrix}
a & & b & & c & & d  \\
 &  \color{red}x&  & \color{red}y &  & \color{red}z &   \\
e & & f & & g & & h \\
 &  \color{red}w&  & \color{red}u &  & \color{red}v &   \\
i & & j & & k & & \ell  \\
 &  \color{red}t&  & \color{red}s &  & \color{red}r &   \\
m & & n & & o & & p 
\end{matrix} \right )
\qquad
\left ( \begin{matrix}
a^* & & b^* & & c^* & & d^*  \\
 &  \color{red}x^*&  & \color{red}y^* &  & \color{red}z^* &   \\
e^* & & f^* & & g^* & & h^* \\
 &  \color{red}w^*&  & \color{red}u^* &  & \color{red}v^* &   \\
i^* & & j^* & & k^* & & \ell^*  \\
 &  \color{red}t^*&  & \color{red}s^* &  & \color{red}r^* &   \\
m^* & & n^* & & o^* & & p^* 
\end{matrix} \right )
\]

The elements $a$,$b$,$c$, etc... belong to the left cumulant matrix $\overline{B}$ while the elements
$x$,$y$,$z$ etc... belong to the left-interlacing matrix $\overline{A}_{\max}$.

Similarly the elements $a^*$,$b^*$,$c^*$, etc... belong to the right cumulant matrix $\underline{B}^*$ while the elements
$x^*$,$y^*$,$z^*$ etc... belong to the right-interlacing matrix $\underline{A}^*_{\min}$.

We wish to show that the value of the entry of $A_{\max}$ at position $u$ is equal to the value of $A^*_{\min}$ at position
$u^*$. That is, by equations (\ref{seconddiff}) and (\ref{seconddiff2}) we want to show that:

\[ x + u - w - y = u^* + z^* - y^* - v^* \]

As a consequence of lemma \ref{lemma} there is some $\gamma$ such that:
\begin{align*}
a + b^* & = b + c^* = c + d^* = \gamma \\
e + f^* & = f + g^* = g + h^* = \gamma+1 \\
i + j^* & = j + k^* = k + \ell^* = \gamma+2 \\
m + n^* & = n + o^* = o + p^* = \gamma+3 \\
\end{align*}

Now, by construction, we have:
\begin{align*}
& x + u - w - y \\
& = \min(a,f-1) + \min(f,k-1) - \min(e,j-1) - \min(b,g-1) \\
& = \min(\gamma - b^*,\gamma - g^*) + \min(\gamma + 1 - g^*, \gamma + 1 - \ell^*) \\
& \phantom{*******} - \min(\gamma + 1 - f^*,\gamma+1- k^*) - \min(\gamma - c^*,\gamma - h^*) \\
& = - \max(b^*,g^*) - \max( g^*, \ell^*) + \max( f^*, k^*) + \max( c^*, h^*) \\
& = -y^* - v^* + u^* + z^*
\end{align*}

The result follows.
\end{proof}
\end{proposition}

A similar argument to the above may be used to show that $A^*_{\max} = A_{\min}$ as well as
$C^*_{\min} = C_{\max}$ and $C^*_{\max} = C_{\min}$.
More precisely:

\begin{proposition}\label{dual}
If $s$ is a binary string, and $\overline{s}$ is its complement, then $A_s = A^*_{\overline{s}}$
and $C_s = C^*_{\overline{s}}$.

\end{proposition}

In other words, for any alternating sign matrix $B$ the partial order of matrices occuring in $\mathfrak{D}^*(B)$ is precisely the dual of the partial order of matrices occuring in $\mathfrak{D}(B)$. Similarly for $\mathfrak{U}(B)$ and $\mathfrak{U}^*(B)$.

\section{Proof of main theorem}\label{proof}

Our proof is almost identical to that given in \cite{robrum}. Let us recall the recurrence:

\begin{align}\label{rec3}
x[k+1]_{i,j} =  \frac{ \mu_{i,n-k-j} x[k]_{i,j} x[k]_{i+1,j+1} + \lambda_{i,j} x[k]_{i,j+1} x[k]_{i+1,j}}{x[k-1]_{i+1,j+1}}
\end{align}

To simplify things, let us introduce the notation:
\[D(x[k])_{i,j} = \mu_{i,n-k-j} x[k]_{ij} x[k]_{i+1,j+1} + \lambda_{ij} \, x[k]_{i,j+1} x[k]_{i+1,j} \]
so that we may rewrite equation \ref{rec3} as:
\[ x[k+1]_{i,j} = \frac{D(x[k])_{i,j}}{s(x[k-1])_{i,j}} \]

\begin{theorem}
For $2 \leq k \leq n$ we have:
\begin{equation} \label{double}
 x[k]_{1,1} = \sum_{\substack{(A,B) \\ |B| = k, |A| = k-1}} 
\lambda^{F(B)} s(\lambda)^{-F(A)} \mu^{G(B)} t(\mu)^{-G(A)}x[1]^B s(x[0])^{-A} 
\end{equation}
The sum is over all pairs of matrix $(A,B)$ such that $A$ occurs in the expansion of $\mathfrak{D}(B)$.

\begin{proof}


The result is trivially true when $k=2$. Making use of the invariance in $k$, followed by the recurrence, we may obtain $x[k+1]_{1,1}$ from $x[k]_{i,j}$ as follows:
\begin{align*} 
x[k+1]_{1,1} & = \sum_{\substack{(A,B) \\ |B| = k, |A| = k-1}} \lambda^{F(B)} s(\lambda)^{-F(A)} \mu^{G(B)} t(\mu)^{-G(A)}x[2]^B s(x[1])^{-A} \\
& = \sum_{\substack{(A,B) \\ |B| = k, |A| = k-1}}  \lambda^{F(B)} s(\lambda)^{-F(A)} \mu^{G(B)} t(\mu)^{-G(A)} \left ( \frac{D(x[1])}{s(x[0])} \right )^B s(x[1])^{-A} 
\end{align*}
We must show that this is equal to:
\[ \sum_{\substack{(B,C) \\ |C| = k+1, |B| = k}} \lambda^{F(C)} s(\lambda)^{-F(B)}\mu^{G(C)} t(\mu)^{-G(B)} x[1]^C s(x[0])^{-B} \]

To do this, we fix some alternating sign matrix $B$ with $|B| = k$ and take the coefficient of $s(x[0])^{-B}$ on both sides. We must now prove that:
\begin{align}
\sum_{ |A| = k-1} & \lambda^{F(B)} s(\lambda)^{-F(A)}  \mu^{G(B)} t(\mu)^{-G(A)} D(x[1])^B s(x[1])^{-A} \label{toprove} \\ &= 
\sum_{|C| = k+1} \lambda^{F(C)} s(\lambda)^{-F(B)} \mu^{G(C)} t(\mu)^{-G(B)} x[1]^C \notag 
\end{align}
Here the sum is over all $A$ (resp. $C$) which may be found in the expansion of $\mathfrak{D}(B)$ (resp $\mathfrak{U}(B)$). 

Making use of proposition \ref{dual} we may rewrite the right hand side of equation (\ref{toprove}) as:
\begin{align}
 \sum_{|C| = k+1} & \lambda^{F(C)} s(\lambda)^{-F(B)} \mu^{G(C)} t(\mu)^{-G(B)} x[1]^C \notag\\ 
& = s(\lambda)^{-F(B)} t(\mu)^{-G(B)} x[1]^{C_{\max}} \lambda^{F(C_{\max})}\mu^{G(C^*_{\min})} 
\prod_{B_{ij} = 1}(\mu_{i,n-k-j} + \lambda_{ij} \frac{x[1]_{i+1,j} x[1]_{i,j+1}}{x[1]_{ij} x[1]_{i+1,j+1}}) \notag\\
 & = s(\lambda)^{-F(B)} t(\mu)^{-G(B)} x[1]^{C_{\max}} \lambda^{F(C_{\max})} \mu^{G(C^*_{\min})} \prod_{B_{ij} = 1} 
\frac{D(x[1]_{i,j})}{x[1]_{ij} s(x[1]_{i,j})}\label{lhs}
\end{align}
while the left hand side of equation (\ref{toprove}) may be written as:
\begin{align}
& \phantom{=} \sum_{ |A| = k-1} \lambda^{F(B)} s(\lambda)^{-F(A)} \mu^{G(B)} t(\mu)^{-G(A)} D(x[1])^B s(x[1])^{-A} \notag\\
& =\lambda^{F(B)} \mu^{G(B)} D(x[1])^B s(\lambda)^{-F(A_{\min})} t(\mu)^{-G(A^*_{\max})} s(x[1])^{-A_{\min}} 
\notag \\ & \phantom{======================}\prod_{B_{i,j} = -1} (\mu_{i,n-k-j} + \lambda_{ij} \frac{x[1]_{i+1,j} x[1]_{i,j+1}}{x[1]_{ij} x[1]_{i+1,j+1}}) \notag\\
& = \lambda^{F(B)} \mu^{G(B)} D(x[1])^B s(\lambda)^{-F(A_{\min})} t(\mu)^{-G(A^*_{\max})} s(x[1])^{-A_{\min}} 
\prod_{B_{i,j} = -1} \frac{D(x[1]_{i,j})}{x[1]_{ij} s(x[1]_{i,j})} \notag \\ \notag
& =  \lambda^{F(B)} \mu^{G(B)} s(\lambda)^{-F(A_{\min})} t(\mu)^{-G(A^*_{\max})} s(x[1])^{-A_{\min}} 
\prod_{B_{i,j} = 1}D(x[1]_{i,j}) \prod_{B_{i,j} = -1} \frac{1}{x[1]_{ij} s(x[1]_{i,j})} \label{rhs}\\
\end{align}

Comparing equation (\ref{lhs}) with equation (\ref{rhs}), we must show that:
\begin{align*} 
s(\lambda)^{F(A_{\min})} & \lambda^{F(C_{\max})}  t(\mu)^{G(A^*_{\max})} \mu^{G(C^*_{\min})} x[1]^{C_{\max}} s(x[1])^{A_{\min}} \\
= & s(\lambda)^{F(B)} \lambda^{F(B)} t(\mu)^{G(B)} \mu^{G(B)} (x[1] s(x[1]))^B
\end{align*}

To complete the proof one need only observe that:
\[ \min(x+1, y) + \max(x,y-1) = x+y \]

More precisely, we have, by construction, that:
\begin{align*}
\overline{A}^{\min}_{i,j} & = \max(\overline{B}_{i,j}, \overline{B}_{i+1,j+1} - 1) \\
\overline{C}^{\max}_{i,j} & = \min(\overline{B}_{i-1,j-1} + 1, \overline{B}_{i,j})
\end{align*}
and so:
\begin{align*} 
\overline{C}^{\max}_{i,j} + \overline{A}^{\min}_{i-1,j-1} 
& = \min(\overline{B}_{i-1,j-1} + 1, \overline{B}_{i,j}) + \max(\overline{B}_{i-1,j-1}, \overline{B}_{i,j} - 1) \\
& = \overline{B}_{i,j} + \overline{B}_{i-1,j-1}
\end{align*}
This gives us the same power of $\lambda_{i,j}$ on both sides. By equations \ref{seconddiff} and \ref{seconddiff2} we also have:
\[ C^{\max}_{i,j} + A^{\min}_{i-1,j-1} = B_{i,j} + B_{i-1,j-1} \]
This gives us the same power of $x_{i,j}$ on both sides.

The power of $\mu_{i,j}$ on the left hand side is given by:
\begin{align*} \underline{C}^{*\min}_{i,k-j+1} + \underline{A}^{*\max}_{i-1,k-j}
& = \max(\underline{B}_{i,k-j+1}, \underline{B}_{i+1,k-j+2}) + 
\min(\underline{B}_{i,k-j+1}, \underline{B}_{i+1,k-j+2}) \\
& = \underline{B}_{i,k-j+1} + \underline{B}_{i+1,k-j+2}
\end{align*}
which is the power of $\mu_{i,j}$ on the right hand side. The result follows.
\end{proof}
\end{theorem}

\nocite{robrum}
\bibliographystyle{alpha}
\bibliography{references}

\end{document}